\documentclass[12pt,A4paper]{article}
\usepackage{latexsym,amssymb,amsmath,amscd,amsfonts,amsthm,mathrsfs}
\usepackage{epsfig,graphics,pstricks,ifpdf}
\usepackage[all]{xy}

\newtheorem{thm}{Theorem}[section]

\newtheorem{prop}[thm]{Proposition}
\newtheorem{lem}[thm]{Lemma}

\def\pf{\noindent{\it Proof.} }
\newcommand{\Integer}{\mathbb{Z}}

\newcommand{\p}{\mathcal{P}}
\def\qed{\nopagebreak\hfill{\rule{4pt}{7pt}}\medbreak}

\makeatletter \@addtoreset{equation}{section} \makeatother

\begin{document}

\begin{center}
{\Large\bf Partitions of $\mathbb{Z}_n$ into Arithmetic
Progressions}
\end{center}

\vskip 2mm

\centerline{William Y.C. Chen$^1$, David G.L. Wang$^2$, and Iris F.
Zhang$^3$}

\begin{center}
Center for Combinatorics, LPMC-TJKLC\\
Nankai University, Tianjin 300071, P.R. China

\vskip 2mm%
$^1$chen@nankai.edu.cn,  $^2$wgl@cfc.nankai.edu.cn,
$^3$zhangfan03@mail.nankai.edu.cn
\end{center}

\begin{abstract}

We introduce the notion of arithmetic progression blocks or
AP-blocks of $\Integer_n$, which can be represented as sequences of
the form $(x, x+m, x+2m, \ldots, x+(i-1)m) \pmod n$. Then we
consider the problem of partitioning $\Integer_n$ into AP-blocks for
a given difference $m$. We show that subject to a technical
condition, the number of partitions of $\Integer_n$ into
$m$-AP-blocks of a given type is independent of $m$.  When we
restrict our attention to blocks of sizes one or two, we are led to
a combinatorial interpretation of a formula recently derived by
Mansour and Sun as a generalization of the Kaplansky numbers. These
numbers have also occurred as the coefficients in Waring's formula
for symmetric functions.

\end{abstract}



\section{Introduction}

Let $\Integer_n$ be the cyclic group of order $n$ whose elements are
written as $1,2,\ldots,n$. Intuitively, we assume that the elements
$1,2,\ldots,n$ are placed clockwise on a cycle. Thus $\Integer_n$
can be viewed as an {\em $n$-cycle}, more specifically, a directed
cycle. In his study of the m\'enages problem, Kaplansky \cite{Kap43}
has shown that the number of ways of choosing $k$ elements from
$\Integer_n$ such that no two elements differ by one modulo $n$ (see
also Brauldi \cite{Bru77}, Comtet \cite{Com74}, Riordan
\cite{Rio58}, Ryser \cite{Rys63} and Stanley \cite[Lemma
2.3.4]{Sta97}) equals
\begin{equation}\label{1}
{n\over n-k}{n-k\choose k}.
\end{equation}
Moreover, Kaplansky \cite{Kap95} considered the following
generalization. Assume that $n\ge pk+1$. Then the number of
$k$-subsets $\{x_1,x_2,\ldots,x_k\}$ of $\Integer_n$ such that
\begin{equation}\label{xij}
x_i-x_j\not\in\{1,2,\ldots,p\}
\end{equation}
for any pair $(x_i, x_j)$ of distinct elements, is given by
\begin{equation}\label{4}
\frac{n}{n-pk}{n-pk\choose k}.
\end{equation}
Here we clarify the meaning of the notation \eqref{xij}. Given two
elements $x$ and $y$ of $\Integer_n$, $x-y$ may be considered as the
distance from $y$ to $x$ on the directed cycle $\Integer_n$.
Therefore, \eqref{xij} says that the distance from any element $x_i$
to any other element $x_j$ on the directed cycle $\Integer_n$ is at
least $p+1$.

From a different perspective, Konvalina \cite{Kon81} studied the
number of $k$-subsets $\{x_1, x_2, \ldots, x_k\}$ such that no two
elements $x_i$ and $x_j$ are ``uni-separated'', namely $x_i-x_j\neq
2$ for all $x_i$ and $x_j$. Remarkably, Konvalina discovered that
the answer is also given by the Kaplansky number \eqref{1} for
$n\ge2k+1$. Other generalizations and related questions have been
investigated by Hwang \cite{Hwa81}, Hwang, Korner and Wei
\cite{Hwa-Kor-Wei84}, Munarini and Salvi \cite{Mun-Sal03}, Prodinger
\cite{Pro83} and Kirschenhofer and Prodinger \cite{Kir86}. Recently,
Mansour and Sun \cite{Man-Sun07} obtained the following unification
of the formulas of Kaplansky and Konvalina.

\begin{thm}\label{Man-Sun} Assume that
$m,p,k\ge1$ and $n\ge mpk+1$. The number of $k$-subsets $\{x_1, x_2,
\ldots, x_k\}$ of $\Integer_n$ such that
\begin{equation}\label{2}
x_i-x_j\not\in\{m,2m,\ldots,pm\}
\end{equation}
for any pair $(x_i,x_j)$, is given by the formula \eqref{4}, and is
independent of $m$.
\end{thm}

In the spirit of the original approach of Kaplansky, Mansour and Sun
first solved the enumeration problem of choosing $k$-subset from an
$n$-set with elements lying on a line. They established a recurrence
relation, and solved the equation by computing the residues of some
Laurent series. The case for an $n$-cycle can be reduced to the case
for a line. They raised the question of finding a combinatorial
proof of their formula. Guo \cite{Guo07} found a proof by using
number theoretic properties and Rothe's identity:
\[
\sum_{k=0}^n\frac{xy}{(x+kz)(y+(n-k)z)}%
{x+kz\choose k}{y+(n-k)z\choose n-k}%
=\frac{x+y}{x+y+nz}{x+y+nz\choose n}.
\]

This paper is motivated by the question of Mansour and Sun. We
introduce the notion of arithmetic progression blocks or AP-blocks
of $\Integer_n$. A sequence of the form
\[
(x, x+m, x+2m, \ldots, x+(i-1)m) \pmod n
\]
is called an AP-block, or an $m$-AP-block, of length $i$
and of difference $m$. Then we consider partitions of
$\Integer_n$ into $m$-AP-blocks $B_1, B_2, \ldots, B_k$ of the same
difference $m$. The type of such a partition is referred to as the
type of the multisets of the sizes of the blocks. Our main result
shows that subject to a technical condition,  the number of
partitions of $\Integer_n$ into $m$-AP-blocks of a given type is
independent of $m$ and is equal to the multinomial coefficient.

This paper is organized as follows. In Section 2, we give a review
of the cycle dissections and make a connection between the Kaplansky
numbers and the cyclic multinomial coefficients. We present the main
result in Section 3, that is, subject to a technical condition, the
number of partitions of $\Integer_n$ into $m$-AP-blocks of a given
type equals the multinomial coefficient and does not depend on $m$.
We present a separation algorithm which leads to a bijection between
$m$-AP-partitions and $m'$-AP-partitions of $\Integer_n$. The
correspondence between $m$-AP-partitions and cycle dissections
($m'=1$) implies the main result Theorem \ref{Main}. For the type
$1^{n-(p+1)k}(p+1)^k$ we are led to a combinatorial proof which
answers the question of  Mansour and Sun.

\section{Cycle Dissections}

\begin{figure}[h]
\begin{center}
\begin{pspicture}(5,2.5)(0,0)
\psset{unit=5pt}%
\pscircle[linewidth=1pt](0,0){10}%
\psarc[linewidth=.5pt]{<-}(0,0){3.5}{0}{120}%
\psarc[linewidth=.5pt]{<-}(0,0){3.5}{120}{240}%
\psarc[linewidth=.5pt]{<-}(0,0){3.5}{240}{360}%
\pscircle*(0,10){.4}%
\pscircle*(3.09,9.51){.4}%
\pscircle*(5.88,8.09){.4}%
\pscircle*(8.09,5.88){.4}%
\pscircle*(9.51,3.09){.4}%
\pscircle*(10,0){.4}%
\pscircle*(9.51,-3.09){.4}%
\pscircle*(8.09,-5.88){.4}%
\pscircle*(5.88,-8.09){.4}%
\pscircle*(3.09,-9.51){.4}%
\pscircle*(0,-10){.4}%
\pscircle*(-3.09,9.51){.4}%
\pscircle*(-5.88,8.09){.4}%
\pscircle*(-8.09,5.88){.4}%
\pscircle*(-9.51,3.09){.4}%
\pscircle*(-10,0){.4}%
\pscircle*(-9.51,-3.09){.4}%
\pscircle*(-8.09,-5.88){.4}%
\pscircle*(-5.88,-8.09){.4}%
\pscircle*(-3.09,-9.51){.4}%
\put(-.5,11.5){\small$1$}%
\put(3,11){\small$2$}%
\put(7,9){\small$3$}%
\put(9.5,6.5){\small$4$}%
\put(10.5,3.5){\small$5$}%
\put(11,-.5){\small$6$}%
 \put(10.5,-4){\small$7$}%
\put(9.5,-7.5){\small$8$}%
\put(7,-10){\small$9$}%
\put(3,-12){\small$10$}%
\put(-1,-13){\small$11$}%
\put(-5,-12.5){\small$12$}%
\put(-9,-10.5){\small$13$}%
\put(-11.5,-8){\small$14$}%
\put(-13,-4){\small$15$}%
\put(-13.5,-.5){\small$16$}%
\put(-13,3){\small$17$}%
\put(-11.5,6.5){\small$18$}%
\put(-9,9.5){\small$19$}%
\put(-5,11){\small$20$}%
\pscircle[linewidth=1pt](30,0){10}%
\psarc[linewidth=.5pt]{<-}(30,0){3.5}{0}{120}%
\psarc[linewidth=.5pt]{<-}(30,0){3.5}{120}{240}%
\psarc[linewidth=.5pt]{<-}(30,0){3.5}{240}{360}%
\pscircle*(30,10){.4}%
\pscircle*(33.09,9.51){.4}%
\pscircle*(35.88,8.09){.4}%
\pscircle*(38.09,5.88){.4}%
\pscircle*(39.51,3.09){.4}%
\pscircle*(40,0){.4}%
\pscircle*(39.51,-3.09){.4}%
\pscircle*(38.09,-5.88){.4}%
\pscircle*(35.88,-8.09){.4}%
\pscircle*(33.09,-9.51){.4}%
\pscircle*(30,-10){.4}%
\pscircle*(26.91,9.51){.4}%
\pscircle*(24.12,8.09){.4}%
\pscircle*(21.91,5.88){.4}%
\pscircle*(20.49,3.09){.4}%
\pscircle*(20,0){.4}%
\pscircle*(20.49,-3.09){.4}%
\pscircle*(21.91,-5.88){.4}%
\pscircle*(24.12,-8.09){.4}%
\pscircle*(26.91,-9.51){.4}%
\psline[linewidth=.5pt](37.574,3.859)(40.247,5.221)%
\psline[linewidth=.5pt](36.01,6.01)(38.132,8.132)%
\psline[linewidth=.5pt](33.859,7.574)(35.221,10.247)%
\psline[linewidth=.5pt](38.395,-1.33)(41.358,-1.799)%
\psline[linewidth=.5pt](37.574,-3.859)(40.247,-5.221)%
\psline[linewidth=.5pt](36.01,-6.01)(38.132,-8.132)%
\psline[linewidth=.5pt](31.33,-8.395)(31.799,-11.358)%
\psline[linewidth=.5pt](22.426,3.859)(19.753,5.221)%
\psline[linewidth=.5pt](23.99,6.01)(21.868,8.132)%
\psline[linewidth=.5pt](26.141,7.574)(24.779,10.247)%
\psline[linewidth=.5pt](21.605,-1.33)(18.642,-1.799)%
\psline[linewidth=.5pt](22.426,-3.859)(19.753,-5.221)%
\psline[linewidth=.5pt](28.67,-8.395)(28.201,-11.358)%
\put(29.5,11.5){\small$1$}%
\put(33,11){\small$2$}%
\put(37,8.89){\small$3$}%
\put(39,6.5){\small$4$}%
\put(40.5,3.5){\small$5$}%
\put(41,-.5){\small$6$}%
\put(40.5,-4){\small$7$}%
\put(39,-7.5){\small$8$}%
\put(37,-10){\small$9$}%
\put(33,-12){\small$10$}%
\put(29,-13){\small$11$}%
\put(25.5,-12.5){\small$12$}%
\put(22,-10.5){\small$13$}%
\put(19,-8){\small$14$}%
\put(17,-4){\small$15$}%
\put(16.5,-.5){\small$16$}%
\put(17,3){\small$17$}%
\put(19,6.5){\small$18$}%
\put(22,9.5){\small$19$}%
\put(25.5,11){\small$20$}%
\end{pspicture}
\end{center}\vspace{50pt}
\caption{A $20$-cycle dissection of type $1^82^33^2$.\label{Fig1}}
\end{figure}
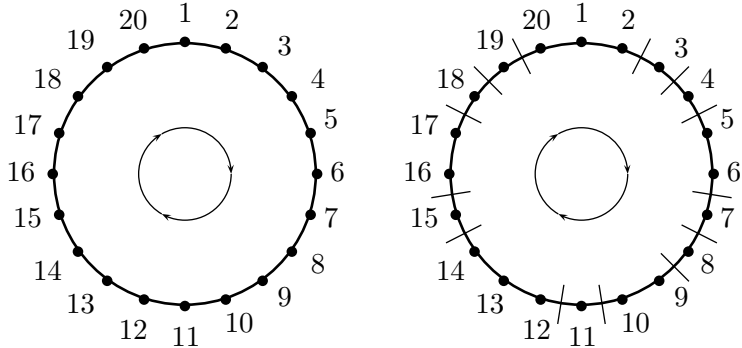

In their combinatorial study of Waring's formula on symmetric
functions, Chen, Lih and Yeh \cite{Chen-Lih-Yeh95} introduced the
notion of cycle dissections. Recall that a {\em dissection of an
$n$-cycle} is a partition of the cycle into blocks, which can be
viewed by putting cutting bars on some edges of the cycle. Note that
there at least one bar to cut a cycle into straight segments. A
dissection of an $n$-cycle is said of {\em type}
$1^{k_1}2^{k_2}\cdots n^{k_n}$ if there are $k_i$ blocks of $i$
elements in it. For instance, Figure \ref{Fig1} gives a $20$-cycle
dissection of type $1^82^33^2$.

The following lemma is due to Chen-Lih-Yeh \cite[Lemma
3.1]{Chen-Lih-Yeh95}.

\begin{lem} \label{lem1}%
For an $n$-cycle, the number of dissections of type
$1^{k_1}2^{k_2}\cdots n^{k_n}$ is given by the {\em cyclic
multinomial coefficients}:
\begin{equation}
{n\over k_1+\cdots+k_n}{k_1+\cdots+k_n\choose k_1,\ldots,k_n}.
\end{equation}
\end{lem}

This lemma is easy to prove. Given a dissection, one may pick up any
segment as a distinguished segment. This can be done in
$k_1+k_2+\cdots+k_n$ ways. On the other hand, any of the $n$
elements can serve as the first element of the distinguished
segment.

Consider a cycle dissection of type $1^{n-(p+1)k}(p+1)^k$. The set
of the first elements of each segment of length $p+1$ corresponds a
$k$-subset of $\Integer_n$ satisfying \eqref{xij}. Thus the cyclic
multinomial coefficient of type $1^{n-(p+1)k}(p+1)^k$ reduces to
\eqref{4} and particularly the cyclic multinomial coefficient of
type $1^{n-2k}2^k$ reduces to the Kaplansky number \eqref{1}.

\section{Partitions of $\Integer_n$ into Arithmetic Progressions}

In this section, we present the main result of this paper, namely, a
formula for the number of partitions of $\Integer_n$ into
$m$-AP-blocks of a given type. The proof is based on a separation
algorithm to transform an $m$-AP-partition to an $m'$-AP-partition.

We begin with some concepts. First, $\Integer_n$ is considered as a
directed cycle. An arithmetic progression block, or an AP-block of
$\Integer_n$, is defined to  be a sequence of elements of
$\Integer_n$ of the following form
\[
B=(x, x+m, x+2m, \ldots, x+(i-1)m) \pmod n,
\]
where $m$ is called the {\em difference} and $i$ is called the {\em
length} of $B$. An AP-block of difference $m$ is called an
$m$-AP-block. If $B$ contains only one element, then it is called a
{\em singleton}. The first element $x$ is called the {\em head} of
$B$. An $m$-AP-partition, or a partition of $\Integer_n$ into
$m$-AP-blocks, is a set of $m$-AP-blocks of $\Integer_n$ whose
underlying sets form a partition of $\Integer_n$. For example,
\begin{equation}\label{eq-ex-1}
(7,9,11),\ (8),\ (10,12),\ (1),\ (2,4,6),\ (3),\ (5)
\end{equation}
is a 2-AP-partition of $\Integer_{12}$ with four singletons and
three non-singleton heads $7$, $10$ and $2$.

It should be noted that different AP-blocks may correspond to the
same underlying set. For example, $(1,3)$ and $(3,1)$ are regarded
as different AP-blocks of $\Integer_4$, but they have the same
underlying set $\{1, 3\}$. On the other hand, as will be seen in
Proposition \ref{AP-block vs block}, it often happens that an
AP-block is uniquely determined by its underlying set. For example,
given the difference $m=3$, the AP-block $(12,15,2,5,8)$ of
$\Integer_{16}$ is uniquely determined by the underlying set
$\{2,5,8,12,15\}$ since there is only one way to order these five
elements to form an arithmetic progression of difference $3$ modulo
16.

For an $m$-AP-partition $\pi$, the {\em type} of $\pi$ is defined by
the type of the multisets of the sizes of the blocks. Usually, we
use the notation $1^{k_1}2^{k_2}\cdots n^{k_n}$ to denote a type for
which there are $k_1$ blocks of size one, $k_2$ blocks of size two,
etc. However, for the sake of presentation, we find it more
convenient to ignore the zero exponents and express a type in the
form $i_1^{k_1}i_2^{k_2}\cdots i_r^{k_r}$, where $1\le
i_1<i_2<\cdots<i_r$ and all $k_j\ge1$. For example, the AP-partition
\eqref{eq-ex-1} is of type $1^42^13^2$.

Throughout this paper, we restrict our attention to
$m$-AP-partitions with at least one singleton block and also at
least one non-singleton block, namely, $i_1=1$ and $r\ge2$ in the
above notation of types. Here is the aforementioned condition:
\begin{equation}\label{TechCondition}
\left\lceil\frac{k_1}{k_2+\cdots+k_r}\right\rceil\ge (m-1)(i_r-1),
\end{equation}
where the notation $\lceil x\rceil$ for a real number $x$
stands for the smallest integer that larger than or
equal to $x$. Obviously, the condition (\ref{TechCondition}) holds
for $m=1$. For $m\ge2$, \eqref{TechCondition} is equivalent to the
relation
\begin{equation}\label{BoundCondition1}
k_1\geq(k_2+\cdots+ k_r)\big[(m-1)(i_r -1)-1\big]+1.
\end{equation}
We prefer the form \eqref{TechCondition} for a reason that will
become clear in the combinatorial argument in the proof of Theorem
\ref{Main}. In fact on an $n$-cycle dissection, the
$\sum_{j=2}^rk_j$ non-singleton heads divide the $k_1$ singletons
into $\sum_{j=2}^rk_j$ segments. By virtue of the pigeonhole
principle, there exists a segment containing at least $(m-1)(i_r-1)$
singletons.

\begin{figure}[h]
\begin{center}
\begin{pspicture}(6,2.5)(0,0)
\psset{unit=5pt} \pscircle[linewidth=1pt](0,0){10}%
\psarc[linewidth=.5pt]{<-}(0,0){3.5}{0}{120}%
\psarc[linewidth=.5pt]{<-}(0,0){3.5}{120}{240}%
\psarc[linewidth=.5pt]{<-}(0,0){3.5}{240}{360}%
\put(-.4,-.4){\small$\pi$}%
\put(6,9.5){\small$1$}%
\put(10,6){\small$2$}%
\put(11.5,0){\small$3$}%
\put(10.5,-6){\small$4$}%
\put(6,-10.5){\small$5$}%
\put(0,-12.5){\small$6$}%
\put(-7,-11){\small$7$}%
\put(-11,-6.5){\small$8$}
\put(-12.5,-.5){\small$9$}%
\put(-11.5,5.5){\small$10$}%
\put(-6,10.5){\small$11$}
\put(0,11.5){\small$12$}%
\pscircle*(5.0,8.5){.3}
\pscircle*(5,-8.5){.3}
\pscircle*(0,-10){.3}
\pscircle*(-8.5,5){.3}
\pscircle[linewidth=.5pt](-5,-8.5){.8}
\pscircle[linewidth=.5pt](-8.5,-5){.8}
\pscircle[linewidth=.5pt](-9.9
,0){.8}
\put(-5.6,8){$\circledast$}
\put(0,9){$\circledast$}
\put(7.5,5){$\circledcirc$}
\put(9,0){$\circledcirc$}
\put(8,-5){$\circledcirc$}
\put(-7,-13){$\blacktriangle$}%
\put(-13,-15.5){\small{\rm starting point}}%
\put(-25,-19.5)%
{\small$(7,8,9),\,(10),\,(11,12),\,(1),\,(2,3,4),\,(5),\,(6)$}%
\psline[linewidth=1pt]{->}(14,0)(20,0)
\psline[linewidth=1pt]{<-}(14,-2)(20,-2)%
\put(16,1){\small$\psi$}%
\put(16,-4){\small$\varphi$}%
\pscircle[linewidth=1pt](35,0){10}
\psarc[linewidth=.5pt]{<-}(35,0){3.5}{0}{120}%
\psarc[linewidth=.5pt]{<-}(35,0){3.5}{120}{240}%
\psarc[linewidth=.5pt]{<-}(35,0){3.5}{240}{360}%
\put(34.5,-.4){\small$\pi'$}%
\put(40.5,9.5){\small$1$}%
\put(44.5,6.5){\small$2$}%
\put(46,-.5){\small$3$}%
\put(45,-6){\small$4$}%
\put(40.5,-10.5){\small$5$}%
\put(34.5,-13){\small$6$}%
\put(28,-11){\small$7$}%
\put(24,-6){\small$8$}%
\put(22,-.4){\small$9$}%
\put(23.5,6){\small$10$}%
\put(28.5,10){\small$11$}%
\put(34,11.5){\small$12$}%
\pscircle*(40,8.5){.3}
\pscircle*(45,0){.3}
\pscircle*(40,-8.5){.3}
\pscircle*(26.5,-5){.3}
\pscircle[linewidth=.5pt](30,-8.5){.8}
\pscircle[linewidth=.5pt](25.1,0){.8}
\pscircle[linewidth=.5pt](30,8.5){.8}
\put(25.7,4.9){$\circledast$}
\put(34,9.4){$\circledast$}
\put(42,5.5){$\circledcirc$}
\put(43,-5){$\circledcirc$}
\put(34,-10.5){$\circledcirc$}
\put(28.5,-13){$\blacktriangle$}%
\put(24,-15.5){\small{\rm starting point}}%
\put(20,-19.5)%
{\small$(7,9,11),\,(8),\,(10,12),\,(1),\,(2,4,6),\,(3),\,(5)$}%
\end{pspicture}
\end{center}\vspace{90pt}
\caption{The algorithms $\psi$ and $\varphi$ for $T=1^42^13^2$,
$m=1$ and $m'=2$.\label{Fig2}}
\end{figure}
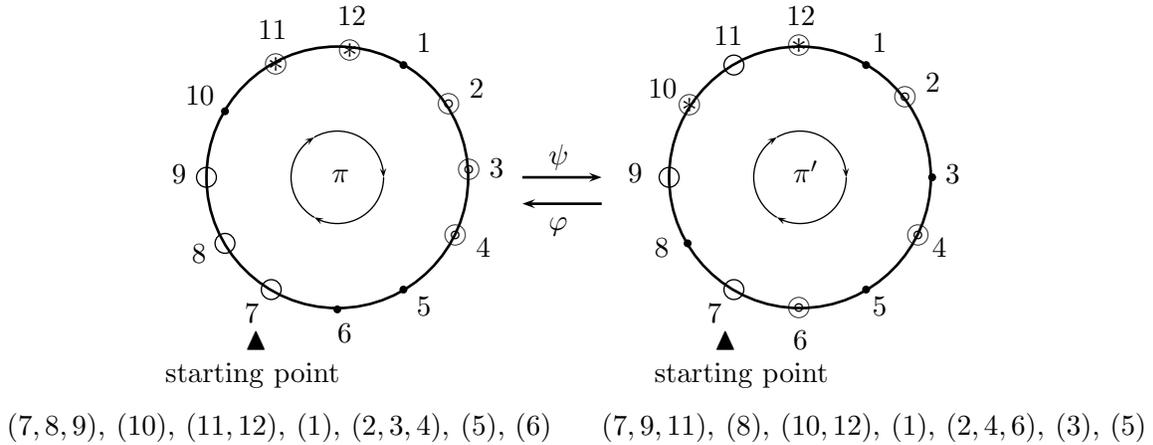

For example in the AP-partition \eqref{eq-ex-1}, the three
non-singleton heads divide the four singletons into three segments
and therefore there exists one segment containing at least $2$
singletons. In this particular partition it is the path from  $2$ to
$7$ that contains two singletons $3$ and $5$, see the right cycle in
Figure \ref{Fig2}.

\begin{prop}\label{AP-block vs block}
Under the condition \eqref{TechCondition}, an $m$-AP-block is not
uniquely determined by its underlying set if and only if $n=i_rm$
and it is of length $i_r$.
\end{prop}

\pf Let $n=i_rm$. Consider the AP-blocks,
\[
B_j=(x+jm,\ x+(j+1)m,\ \ldots,\ x+(j+i_r-1)m) \pmod n, \quad 0\le
j\le i_r-1.
\]
It is easy to see that these AP-blocks $B_j$ $(j=0, 1,\ldots,
i_r-1)$ have the same underlying set
\[
\{x,\ x+m,\ \ldots,\ x+(i_r-1)m\}.
\]

Conversely, suppose that there is an $m$-AP-block $B$ of length
$i_s$ which is not uniquely determined by its underlying set. We may
assume that there exists another AP-block $B'$ having the same
underlying set as $B$. Thus the difference between $B$ and $B'$ lies
only in the order of their elements as a sequence. It follows that
$n=i_sm$ for some $2\le s\le r$. If $m=1$, then $n=i_s$ which yields
$s=r=1$, a contradiction. So we may assume that $m\ge2$ and $2\le
s\le r-1$. Hence $i_s\le i_{r-1}\le i_r-1$, and so
\[
k_1+\sum_{j=2}^rk_ji_j=n=i_sm\le(i_r-1)m.
\]
In view of the condition \eqref{BoundCondition1}, we deduce
that
\[
(i_r-1)m-\sum_{j=2}^rk_ji_j%
\ge k_1\ge\left[(m-1)(i_r-1)-1\right]\sum_{j=2}^rk_j+1
\]
which can be rewritten as
\[
1+\sum_{j=2}^{r-1}k_ji_j+(i_r-1)m\left(\sum_{j=2}^rk_j-1\right)%
\le i_r\sum_{j=2}^{r-1}k_j.
\]
Clearly,
\[
\sum_{j=2}^rk_j-1\ge\sum_{j=2}^{r-1}k_j,
\]
so $(i_r-1)m<i_r$ and thus $i_r<m/(m-1)\le2$ which implies $i_r=1$,
a contradiction. Thus we conclude that $s=r$. This completes the proof.
\qed

For example, the AP-partition \eqref{eq-ex-1} is uniquely determined
by its underlying partition:
\[
\{7,9,11\},\ \{8\},\ \{10,12\},\ \{1\},\ \{2,4,6\},\ \{3\},\ \{5\}.
\]

We are now ready to present the main result of this paper.

\begin{thm}\label{Main}
Given a type $1^{k_1} i_2^{k_2} \cdots i_r^{k_r}$ satisfying the
condition \eqref{TechCondition}, the number of $m$-AP-partitions of
$\Integer_n$  does not depend on $m$, and is equal to the cyclic
multinomial coefficient
\begin{equation}\label{CMC}
{n\over k_1+\cdots+k_r}{k_1+\cdots+k_r\choose k_1,\ldots,k_r}.
\end{equation}
\end{thm}

In fact, Theorem \ref{Main} reduces to Theorem \ref{Man-Sun} when we
specialize the type to $1^{n-(p+1)k}(p+1)^k$. In this case the
condition \eqref{TechCondition} becomes $n\ge kmp+1$. The heads of
the $k$ AP-blocks of length $p+1$ satisfy the condition \eqref{2}.
Conversely, any $k$-subset of $\Integer_n$ satisfying \eqref{2}
determines an $m$-AP-partition of the given type. The cyclic
multinomial coefficient \eqref{CMC} agrees with the formula
\eqref{4} of Theorem \ref{Man-Sun}. For example,  given the type
$1^42^13^2$ and difference $2$, the AP-partition \eqref{eq-ex-1} is
determined by the selection of $\{7,10,2\}$ as heads from
$\Integer_{12}$.

Note that the cyclic multinomial coefficient \eqref{CMC} has
occurred in Lemma \ref{lem1}. Indeed, Lemma 1 is the special case of
Theorem \ref{Main} for $m=1$. We proceed to describe an algorithm,
called the {\em separation algorithm}, to transform
$m$-AP-partitions to $m'$-AP-partitions of the same type
$T=i_1^{k_1}i_2^{k_2}\cdots i_r^{k_r}$, assuming the following
condition holds:
\begin{equation}\label{TechConditionGeneral}
\left\lceil\frac{k_1}{k_2+\cdots+k_r}\right\rceil\ge
(\max\{m,m'\}-1)(i_r-1).
\end{equation}

The separation algorithm enables us to verify Theorem \ref{Main}. We
will state our algorithm for $m$-AP-partitions and
$m'$-AP-partitions, instead of restricting $m'$ to one, because it
is more convenient to present the proof by exchanging the role of
$m$ and $m'$.

Given a type $T=1^{k_1}i_2^{k_2}\cdots i_r^{k_r}$, let $\p_m$ be the
set of $m$-AP-partitions of type $T$. To prove Theorem \ref{Main},
it suffices to show that there is a bijection between $\p_m$ and
$\p_m'$ under the condition \eqref{TechConditionGeneral}.

Let $\pi\in\p_m$. Denote by $H(\pi)$ the set of heads in $\pi$. For
each head $h$ of $\pi$, we consider the nearest non-singleton head
in the counterclockwise direction, denoted $h^*$. Then we denote by
$g(h)$ the number of singletons lying on the path from $h^*$ to $h$
under the convention that $h$ is not counted by $g(h)$. For example,
for the AP-partition $\pi'$ on the right of Figure \ref{Fig2}, we
have $H(\pi')=\{1,2,3,5,7,8,10\}$, $g(1)=g(3)=g(8)=0$,
$g(2)=g(5)=g(10)=1$ and $g(7)=2$. The values $g(h)$ will be needed
in the separation algorithm.

\noindent {\bf The Separation Algorithm.} Let $\pi$ be an
$m$-AP-partition of type $T$. As the first step, we choose a head
$h_1$ of $\pi$, called the {\em starting point}, such that $g(h_1)$
is the maximum. Then we impose a linear order on the elements of
$\Integer_n$ with respect to the choice of $h_1$:
\begin{equation}\label{order}
h_1<h_1+1<h_1+2<\cdots<h_1-1 \quad \pmod n.
\end{equation}
In accordance with the above order, we denote the heads of $\pi$ by
$h_1< h_2< \cdots < h_{t}$, where $t=\sum_{i=1}^rk_i$. The
$m$-AP-block of $\pi$ with head $h_i$ is denoted by $B_i$. Let $l_i$
be the length of $B_i$, and so $\sum_{i=1}^{t}l_i=n$.

We now aim to construct $m'$-AP-blocks $B_1', B_2', \ldots, B_t'$
such that $B_i'$ has the same number of elements as $B_i$. We begin
with $B_1'$ by setting $h_1'=h_1$ and letting $B_1'$ be the
$m'$-AP-block of length $l_1$, namely,
\[
B_1'=\left(%
h_1',\ h_1'+m',\ \ldots,\ h_1'+(l_1-1)m'
\right).%
\]
Among the remaining elements, namely, those that are not in $B_1'$,
we choose the smallest element with respect to \eqref{order},
denoted by $h_2'$, and let $B_2'$ be the $m'$-AP-block of length
$l_2$ with head $h_2'$. Repeating the above procedure, as will be
justified later, after $t$ steps we obtain an $m'$-AP-partition,
denoted $\psi(\pi)$, of type $T$ with blocks $B_1', B_2', \ldots,
B_t'$.

Figure \ref{Fig2} illustrates the separation algorithm from a
$1$-AP-partition $\pi$ to a $2$-AP-partition $\pi'$ of the same type
$T=1^42^13^2$ and vice versa. The solid dots stand for singletons,
whereas the other symbols represent different AP-blocks.

We remark that, as indicated by the example, the starting point can
never be a singleton. In fact, if $s$ is a singleton and $h$ is a
non-singleton head such that all the heads lying on the path from
$s$ to $h$ are singletons, then we have the relation $g(h)>g(s)$.
Since $g(h_1)$ is  maximum, we see that the starting point is always
a non-singleton head.

Clearly, it is necessary to demonstrate that the above algorithm
$\psi$ is valid, namely, we need to justify that underlying sets of
the blocks $B_1', B_2', \ldots, B_t'$ are disjoint.

\begin{prop}\label{well-defined}
The mapping $\psi$ is well-defined, and for any $\pi\in\p_m$, we
have $\psi(\pi)\in \p_{m'}$.
\end{prop}

\pf Let $\pi\in \p_m$ with AP-blocks $B_1, B_2, \ldots, B_t$.
Without loss of generality, we may assume that $h_1, h_2, \ldots,
h_t$ are the heads of $B_1, B_2, \ldots, B_t$, where $h_1$ is the
starting point for the mapping $\psi$ and $h_1', h_2', \ldots, h_t'$
are the corresponding heads generated by $\psi$. Let $l_i$ be the
length of $B_i$. Suppose to the contrary that there exist two heads
$h_i$ and $h_j$ ($i<j$) such that
\[
h_i'+am'\equiv h_j'+bm' \pmod n,
\]
where $0\le a\le l_i-1$ and $0\le b\le l_j-1$.

If $a\ge b$, then $0\le a-b\le l_i-1$ and $h_j'\equiv
h_i'+(a-b)m'\pmod n$. But the point $h_i'+(a-b)m'$ is in $B_i'$,
contradicting the choice of $h_j'$. This yields $a<b$ and thus $0\le
b-a \le l_j-1$.

We claim that the starting point $h_1$ lies on the path from $h_j'$
to $h_i'$. In fact, when the Algorithm $\psi$ is at the $j$-th step
to deal with the head $h_j$, all the points smaller than $h_i'$ lie
in one of the blocks $B_1', B_2', \ldots, B_i'$. Then we see that
$h_j'>h_i'$. Meanwhile, there are  $n-l_1-l_2-\cdots-l_{j-1}>0$
points which are not contained in $B_1', B_2',\ldots, B_{j-1}'$. But
the head $h_j'$ is chosen to be the smallest point not in
$B_1',B_2', \ldots, B_{j-1}'$, we find that $h_j'$ lies on the path
from $h_i'$ to $h_1$.

In addition to $h_i'$ and $h_j'$, we assume that there are $N$
points on the path from $h_j'$ to $h_i'$. Since $h_i'\equiv
h_j'+(b-a)m'\pmod n$ and $1\le b-a\le l_j-1$, we obtain
$N=(b-a)m'-1$. On the other hand, at the $j$-th step, in addition to
the point $h_j'$, there are at least $l_j-1$ points not contained in
$B_1', B_2', \ldots, B_{j-1}'$. Similarly, the choice of $h_1$ and
the condition \eqref{TechConditionGeneral} yield that the largest
$(\max\{m,m'\}-1)(i_r-1)$ heads with respect to the order
\eqref{order} are all singletons by the pigeonhole principle.
Therefore, there are at least $(\max\{m,m'\}-1)(i_r-1)$ points not
contained in $B_1', B_2', \ldots, B_{j-1}'$. It follows that
\begin{equation}
N\ge(\max\{m,m'\}-1)(i_r-1)+(l_j-1).
\end{equation}
Since $N=(b-a)m'-1$ and $1\le b-a\le l_j-1$, we deduce that
\[
(m'-1)(i_r-1)+(l_j-1)\le(b-a)m'-1\le(l_j-1)m'-1,
\]
leading to the contradiction $l_j>i_r$. This completes the proof.
\qed

\begin{prop}\label{prop1} Given an $m$-AP-partition of
$\Integer_n$, the separation algorithm $\psi$ generates the same
$m'$-AP-partition regardless of the choice of the starting point
subject to the maximum property.
\end{prop}

\pf Let $\pi$ be an $m$-AP-partition of $\Integer_n$. Suppose that
$u_1,u_2,\ldots,u_s$ $(s\ge 2)$ are all the heads such that
$g(u_1)=g(u_2)=\cdots=g(u_s)$ is the maximum on $\pi$. Let $u_1$ be
the starting point and $u_1<u_2<\cdots<u_s$ with respect to
\eqref{order}.

It suffices to show that when the Algorithm $\psi$ processes $u_i$
($1\le i\le s$), the $m'$-AP-blocks which have been generated
consist of  all the elements smaller than $u_i$. By induction we
assume that this statement holds  up to $u_{j-1}$.

Let $v_q,v_{q-1},\ldots,v_1,u_j$ be all heads lying on the path $Q$
from $u_{j-1}$ to $u_j$ such that
$u_{j-1}=v_q<v_{q-1}<\cdots<v_1<u_j$. Let $B_i$ be the $m$-AP-block
containing $v_i$. Let $l_i$ be the length of $B_i$ and
\[
B_i'=(v_i',v_i'+m',\ldots,v_i'+(l_i-1)m')
\]
be the corresponding $m'$-AP-blocks generated by the Algorithm
$\psi$. It suffices to show that the path $Q$ consists of the
elements of $B_s',B_{s-1}',\ldots,B_1'$.

Suppose that $v_1,v_2,\ldots,v_p$ are all singletons, but $v_{p+1}$
is not a singleton. Then $p\le q-1$ since $u_{j-1}$ is always a
non-singleton head. The condition \eqref{TechConditionGeneral}
yields that
\[
p\ge(\max\{m,m'\}-1)(i_r-1).
\]

We now wish to show that for any $1\le i\le q$, the block $B_i$ lies
entirely on the path $Q$. If $i\le p$, then $B_i=(v_i)$ is a
singleton block lying  on $Q$. Otherwise, we have $i\ge p+1$ and
\[
B_i=(v_i,v_i+m,\ldots,v_i+(l_i-1)m).
\]
But the total number of points between any two consecutive elements
of $B_i$ is
\[
(l_i-1)(m-1)\le(\max\{m,m'\}-1)(i_r-1)\le p.
\]
Intuitively, all these points can be fulfilled by the singletons
$v_p,v_{p-1},\ldots,v_1$. Since $u_j>v_1$, the largest element
$v_i+(l_i-1)m$ in the block $B_i$ is smaller than $u_j$. Hence the
block $B_i$ ($i=1,2,\ldots,q$) lies entirely on $Q$.

Therefore,  the total number of elements in $B_q,B_{q-1},\ldots,B_1$
equals the length $u_j-u_{j-1}$ of the path $Q$. Since  $B_i'$ has
the same number of elements as $B_i$, the total number of elements
in $B_q',B_{q-1}',\ldots,B_1'$ also equals $u_j-u_{j-1}$.

Moreover, it can be shown that the block $B_i'$ also lies entirely
on the path $Q$ for any $1\le i\le q$. If $i\le p$, the block
$B_i'=(v_i')$ is a singleton given by the separation algorithm.
Since the total number of elements in
$B_q',B_{q-1}',\ldots,B_{i+1}'$ is smaller than $u_j-u_{j-1}$ and
$v_i'$ is chosen to be the smallest element which is not in
$B_q',B_{q-1}',\ldots,B_{i+1}'$, we see the relation $v_i'<u_j$.
Otherwise, we have $i\ge p+1$ and the total number of points between
any two consecutive elements of $B_i'$ equals
\[
(l_i-1)(m'-1)\le(\max\{m,m'\}-1)(i_r-1)\le p.
\]
Intuitively, all these points can be fulfilled by the singletons
$v_p',v_{p-1}',\ldots,v_1'$. Since $u_j>v_1'$, the largest element
$v_i'+(l_i-1)m'$ in the block $B_i'$ is smaller than $u_j$.
Consequently, the block $B_i'$ lies entirely on $Q$.

In summary,  the total number of elements in
$B_q',B_{q-1}',\ldots,B_1'$ which lie on the path $Q$ coincides with
the length of $Q$. Hence the path $Q$ consists  of the elements of
$B_s',B_{s-1}',\ldots,B_1'$. This completes the proof. \qed

\begin{thm}\label{thm4}
Let $T$ be a type as given before. The separation algorithm induces
a bijection between $\p_m$ and $\p_{m'}$ under the condition
\eqref{TechConditionGeneral}.
\end{thm}

\pf We may employ the separation algorithm by interchanging the
roles of $m$ and $m'$ to construct an $m$-AP-partition from an
$m'$-AP-partition, and we denote this map by $\varphi$. We aim to
show that $\varphi$ is indeed the inverse map of $\psi$, namely,
$\varphi(\psi(\pi))=\pi$ for any $\pi\in\p_m$.

Let $h_1,h_2,\ldots,h_t$ be the heads of $\pi$ for the map $\psi$,
where $h_1$ is the starting point. Assume that $\pi$ has AP-blocks
$B_1, B_2, \ldots, B_t$  with $h_i$ being the head of $B_i$.  Let
$l_i$ be the length of $B_i$. By the construction of $\psi$, the
generated heads $h_1'=h_1, h_2', \ldots, h_t'$ have the order
$h_1'<h_2'<\cdots<h_t'$ in accordance with $h_1<h_2<\cdots<h_t$. It
follows that $g(h_1')$ is the maximum considering all heads of the
AP-partition $\psi(\pi)$.

We now apply the map $\varphi$ on the $m'$-AP-partition $\psi(\pi)$
and choose $h_1'$ as the starting point. Let
$h_1'',h_2'',\ldots,h_t''$ be the heads generated by $\varphi$
respectively. In light of  the construction of $\varphi$, we have
$h_1''=h_1'=h_1$ and $h_1''<h_2''<\cdots<h_t''$.

For any $i$, the separation algorithm has the property that the
length of the $m$-AP-block in $\varphi(\psi(\pi))$ containing
$h_i''$ is $l_i$, which is the length of the $m$-AP-block in $\pi$
containing $h_i$.

Note that both $\varphi(\psi(\pi))$ and $\pi$ are $m$-AP-partitions.
They have the same starting point $h_1''=h_1$ and the same length
sequence $(l_1,l_2,\ldots,l_t)$. Thus for any $i=2,3,\ldots,t$, the
head $h_i''$ is the smallest point which is not contained in the
$m$-AP-blocks $B_1,B_2,\ldots,B_{i-1}$, and so does $h_i$. Hence we
conclude that $h_i''=h_i$ and $\varphi(\psi(\pi))=\pi$. This
completes the proof. \qed

\noindent{\bf Acknowledgments.} This work was supported by the 973
Project, the PCSIRT Project of the Ministry of Education, the
Ministry of Science and Technology, and the National Science
Foundation of China.



\begin{thebibliography}{99}
\small \setlength{\itemsep}{-.8mm}

\bibitem{Bru77}%
R.A. Brualdi, Introductory Combinatorics, North-Holland, New York,
1977.

\bibitem{Chen-Lih-Yeh95}%
W.Y.C. Chen, K.W. Lih and Y.N. Yeh, Cyclic tableaux and symmetric
functions, Studies in Applied Math.  94 (1995) 327--339.

\bibitem{Com74}%
L. Comtet, Advanced Combinatorics, D. Reidel Pub. Co., Dordrecht,
Holland, 1974.

\bibitem{Guo07}%
V.J.W. Guo, A new proof of a theorem of Mansour and Sun, European J.
Combin. (2007), to appear.

\bibitem{Hwa81}%
F.K. Hwang, Cycle polynomials, Proc. Amer. Math. Soc., Vol. 83, No.
1 (1981) 215--219.

\bibitem{Hwa-Kor-Wei84}%
F.K. Hwang, J. Korner, and V.K.-W. Wei, Selecting non-consecutive
balls arranged in many lines, J. Combin. Theory Ser. A  37 (1984)
327--336.

\bibitem{Kap43}%
I. Kaplansky, Solution of the ``Probl\`eme des m\'enages'', Bull.
Amer. Math. Soc.  49 (1943) 784--785.

\bibitem{Kap95}%
I. Kaplansky, Selected Papers and Other Writings, Springer (1995)
25--26.

\bibitem{Kir86}%
P. Kirschenhofer and H. Prodinger, Two selection problems revisited,
J. Combin. Theory Ser. A  42 (1986) 310--316.

\bibitem{Kon81}%
J. Konvalina, On the number of combinations without unit separation,
J. Combin. Theory, Ser. A  31 (1981) 101--107.

\bibitem{Man-Sun07}%
T. Mansour and Y. Sun, On the number of combinatorics without
certain separations, European J. Combin. 29 £¨5£© (2008) 1200-1206.

\bibitem{Mun-Sal03}%
E. Munarini, N.Z. Salvi, Scattered subsets, Discrete Math.  267
(2003) 213--228.

\bibitem{Pro83}%
H. Prodinger, On the number combinations without a fixed distance,
J. Combin. Theory Ser. A  35 (1983) 362--365.

\bibitem{Rio58}%
J. Riordan, An Introduction to Combinatorial Analysis, Wiley, New
York, 1958.

\bibitem{Rys63}%
H.J. Ryser, Combinatorial Mathematics, Carus Monograph  14,
Mathematical Association of America, Wiley, New York, 1963.

\bibitem{Sta97}%
R.P. Stanley, Enumerative Combinatorics, Vol.  1, 2nd ed.,
Cambridge, New York, Cambridge University Press, 1997.

\end{thebibliography}
\end{document}